%% file: MSym_extabs.tex
\documentclass[11pt]{amsart}
\usepackage{amssymb, colordvi,multicol}
\usepackage{graphicx,xspace}
\usepackage{amsmath,amssymb,amsfonts,bbm}
\usepackage{amsbsy}
\usepackage[all]{xy}
\usepackage{pstricks,pst-node,pst-tree,pst-poly,pst-text}
\usepackage[width=.85\textwidth,font=small,labelfont=sc]{caption}
\usepackage[bookmarks=true,colorlinks,linkcolor=mylinkcolor,citecolor=mycitecolor]{hyperref} 
\usepackage{bookmark}
	\definecolor{mycitecolor}{rgb}{.2,.7,.2}
	\definecolor{mylinkcolor}{rgb}{0,0,.8}

\newtheorem{theo}{Theorem}
\newtheorem{prop}[theo]{Proposition}

\newtheorem{coro}[theo]{Corollary}

\newtheorem{definition}[theo]{Definition}
\newtheorem{example}[theo]{Example}
\newtheorem{remark}[theo]{Remark}
\newenvironment{defi}{\rm\begin{definition}\rm}{\end{definition}}
\newenvironment{exam}{\rm\begin{example}\rm}{\end{example}}
\newenvironment{rema}{\rm\begin{remark}\rm}{\end{remark}}

\newcounter{FNC}[page]
\def\fauxfootnote#1{{\addtocounter{FNC}{2}$^\fnsymbol{FNC}$%
     \let\thefootnote\relax\footnotetext{$^\fnsymbol{FNC}$#1}}}

\newpsobject{showgrid}{psgrid}{%
  subgriddiv=1,griddots=10,gridlabels=6pt}
\newif\ifhrule\hrulefalse

\def\N{\mathbb N}

\def\Sym{\mathfrak S}
\def\s{\Sym}
\def\m{\mathcal M}
\def\y{\mathcal Y}

\def\ssym{\mathfrak SSym}
\def\msym{\mathcal MSym}
\def\ysym{\mathcal YSym}
\def\qsym{\mathcal QSym}

\def\subGalois{interval retract\xspace}
\def\id{\mathbbm{1}} 
\def\btau{\bs\tau}
\def\bbeta{\bs\beta}
\def\bphi{\bs\phi}
\def\brho{\bs\rho}

\def\bb{\boldbullet}
	\def\boldbullet{{\!\mbox{\LARGE$\mathbf\cdot$}}}  
\def\psplit{\stackrel{\curlyvee}{\to}}
\def\rsplit{\stackrel{\curlyvee_{\!\!+}}{\longrightarrow}}

\newcommand{\Max}[1][]{\mathsf{max}\if#1\else{\left(#1\right)}\fi}
\newcommand{\Min}[1][]{\mathsf{min}\if#1\else{\left(#1\right)}\fi}
\def\maxmax{\hbox{\small$\mathsf{MM}$}}
\def\maxmin{\hbox{\small$\mathsf{Mm}$}}

\def\minmin{\hbox{\small$\mathsf{mm}$}}

\def\bs{\boldsymbol}
\def\into{\hookrightarrow}
\def\onto{\twoheadrightarrow}

\newcommand{\Hilb}[2]{\mathrm{Hilb}_{#1}(#2)}
\newcommand{\demph}[1]{{\Blue{\emph{#1}}}}

\newcommand{\newatop}[2]{\genfrac{}{}{0pt}{}{#1}{#2}}

\newdir^{ (}{{}*!/-5pt/@^{(}}%

\title[Hopf Structures on Trees]{%
New Hopf Structures on Binary Trees \\(Extended Abstract)
} 

\author{Stefan Forcey} 
\address[Forcey]{
	Department of Mathematics\\
       Tennessee State University\\
       3500 John A Merritt Blvd\\
       Nashville, Tennessee \ 37209 
         }
\email{sforcey@tnstate.edu} 
  \urladdr{http://faculty.tnstate.edu/sforcey/}

\author{Aaron Lauve}
\address[Lauve]{
	Department of Mathematics\\ 
       Texas A\&M University\\
       College Station, TX \ 77843 
         }
\email{lauve@math.tamu.edu}
  \urladdr{http://www.math.tamu.edu/\~{}lauve/}

\author{Frank Sottile}
\address[Sottile]{
	Department of Mathematics\\
        Texas A\&M University\\
        College Station, TX \ 77843 
         }
\email{sottile@math.tamu.edu}
\urladdr{http://www.math.tamu.edu/\~{}sottile/}
\thanks{Sottile supported by the NSF grant DMS-0701050.}  

\date{25 August 2009}

\begin{document} 

\begin{abstract} 
The multiplihedra $\m_{\bb} = (\m_n)_{n\geq1}$ form a family of polytopes originating in the study of higher categories and homotopy theory. While the multiplihedra may be unfamiliar to the algebraic combinatorics community, it is nestled between two families of polytopes that certainly are not: the permutahedra $\s_\bb$ and associahedra $\y_\bb$.
The maps $\s_\bb \onto \m_\bb \onto \y_\bb$ reveal several new Hopf structures on tree-like objects nestled between the Hopf algebras $\ssym$ and $\ysym$. We begin their study here, showing that $\msym$ is a module over $\ssym$ and a Hopf module over $\ysym$. Rich structural information about $\msym$ is uncovered via a change of basis---using M\"obius inversion in posets built on the $1$-skeleta of $\m_\bb$. Our analysis uses the notion of an interval retract, which should have independent interest in poset combinatorics. It also reveals new families of polytopes, and even a new factorization of a known projection from the associahedra to hypercubes. \\[0ex]

\noindent{\sc R\'esum\'e.}\ %
Les multiplih\'edra $\m_{\bb} = (\m_n)_{n\geq1}$ formez une famille des polytopes provenant de l'\'etude des cat\'egories plus \'elev\'ees et de la th\'eorie homotopy. Tandis que le multiplih\'edra peut \^etre peu familier \`a la communaut\'e alg\'ebrique de combinatoire, il est nich\'e entre deux familles des polytopes qui ne sont pas certainement : le permutah\'edra  $\s_\bb$ et l'associah\'edra $\y_\bb$.
Les morphisms $\s_\bb \onto \m_\bb \onto \y_\bb$ indiquent plusieurs nouvelles structures de Hopf en fonction arbre-comme des objets nich\'es entre les alg\`ebres de Hopfs $\ssym$ et $\ysym$. Nous commen\,cons leur \'etude ici, prouvant que $\msym$ est un module au-dessus de $\ssym$ et un module de Hopf au-dessus de $\ysym$. Des informations structurales riches sur $\msym$ sont d\'ecouvertes par l'interm\'ediaire d'une modification de base---utilisant inversion de M\"obius dans les posets \'etablis sur le $1$-skeleta de $\m_\bb$. Notre analyse utilise la notion d'un intervalle se r\'etractent, qui devrait avoir l'int\'er\^et ind\'ependant pour la combinatoire de poset. Elle indique \'egalement de nouvelles familles des polytopes, et m\^eme une nouvelle factorisation d'une projection connue de l'associah\'edra aux hypercubes.
\end{abstract}

\maketitle

\section*{Introduction}\label{sec: intro}

In the past 30 years, there has been an explosion of interest in combinatorial Hopf algebras related to the classical ring of symmetric functions. This is due in part to their applications in combinatorics and representation theory, but also in part to a viewpoint expressed in the elegant commuting diagram
\begin{gather*}
\xymatrix@=4pt@C=6pt@R=3pt{
{NSym} \ar@{->>}[ddd]_{} \ar@{^{ (}->}[rrr]_{} &&& {\ssym}  \ar@{->>}[ddd]_{} \\  
\mbox{ } \\
\mbox{ } \\
{Sym} \ar@{^{ (}->}[rrr]_{} &&& {\qsym} 
}\raisebox{-4ex}{.}
\end{gather*}
Namely, much information about an object may be gained by studying how it interacts with its surroundings. From this picture, we focus on the right edge,  $\ssym \onto \qsym$.
We factor this map through finer and finer structures (some well-known and some new) until this edge is replaced by a veritable zoo of Hopf structures. A surprising feature of our results is that each of these factorizations may be given geometric meaning---they correspond to successive polytope quotients from permutahedra to hypercubes.

\subsection*{The (known) cast of characters} Let us reacquaint ourselves with some of the characters who have already appeared on stage. 
\smallskip

{$\ssym$}\ --\ the Hopf algebra introduced by Malvenuto and Reutenauer \cite{MalReu:1995} to explain the isomorphism $\qsym \simeq ({NSym})^*$. A graded, noncommutative, noncocommutative, self-dual Hopf algebra, with basis indexed by permutations,  
it offers a natural setting to practice noncommutative character theory \cite{BleSch:2005}. 

\smallskip

{$\ysym$}\ --\ the (dual of the) Hopf algebra of trees introduced by Loday and Ronco \cite{LodRon:1998}. A graded, noncommutative, noncocommutative Hopf algebra with basis indexed by planar binary trees, it is important for its connections to the Connes-Kreimer renormalization procedure.\smallskip

{$\qsym$}\ --\ The Hopf algebra of quasisymmetric functions introduced by Gessel \cite{Ges:1984} in his study of $P$-partitions. A graded, commutative, noncocommutative Hopf algebra with basis indexed by compositions, it holds a special place in the world of combinatorial Hopf algebras \cite{AguBerSot:2006}.

\subsection*{The new players} In this extended abstract, we study in detail a family of planar binary trees that we call \emph{bi-leveled trees,} which possess two types of internal nodes (circled or not, subject to certain rules). 
These objects are the vertices of Stasheff's multiplihedra \cite{Sta:1970}, originating from his study of $A_\infty$ categories. The multiplihedra were given the structure of CW-complexes by Iwase and Mimura \cite{IwaMim:1989} and realized as polytopes later \cite{For:2008}. They persist as important objects of study, among other reasons, because they catalog all possible ways to multiply objects in the domain and range of a function $f$, when both have nonassociative multiplication rules. More recently, they have appeared as moduli spaces of ``stable quilted discs'' \cite{MauWoo:xx}.

In Section \ref{sec: msym}, we define a vector space $\msym$ with basis indexed by these bi-leveled trees. We give $\msym$ a module structure for $\ssym$ by virtue of the factorization 
$$\ssym \stackrel{\bbeta}{\relbar\joinrel\onto} \msym \stackrel{\bphi}{\relbar\joinrel\onto} \ysym$$ (evident on the level of planar binary trees) and a splitting $\msym \hookrightarrow \ssym$. We also show that $\msym$ is a Hopf module for $\ysym$ and we give an explicit realization of the fundamental theorem of Hopf modules. That is, we find the coinvariants for this action. Our proof, sketched in Section \ref{sec: proof}, rests on a result about poset maps of independent interest.

We conclude in Section \ref{sec: polytopes} with a massive commuting diagram---containing several new families of planar binary trees---that further factors the map from $\ssym$ to $\qsym$. The remarkable feature of this diagram is that it comes from polytopes (some of them even new) and successive polytope quotients. Careful study of the interplay between the algebra and geometry will be carried out in future work.

\section{Basic combinatorial data}\label{sec: basic}

\subsection{Ordered and planar binary trees}\label{sec: trees}
We recall a map $\tau$ from permutations $\s_\bb=\bigcup_n \s_n$ to planar binary trees $\y_\bb=\bigcup_n \y_n$ that has proven useful in many contexts \cite{Ton:1997,LodRon:2002}. Its behavior is best described in the reverse direction as follows. Fix a tree $t \in \y_n$. The $n$ internal nodes of $t$ are equipped with a partial order, viewing the root node as maximal. An \demph{ordered tree} is a planar binary tree, together with a linear extension of the poset of its nodes. These are in bijection with permutations, as the nodes are naturally indexed left-to-right by the numbers $1,\dotsc,n$. The map $\tau$ takes an ordered tree (permutation) to the unique tree whose partial order it extends.

\begin{exam}\label{ex: tau} 
The permutations $1423$, $2413$, and $3412$ share a common image under $\tau$:
\begin{gather*}
	\tau^{-1}\left(\hskip.1em\raisebox{-3ex}{$\includegraphics{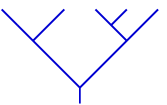}$}\hskip.1em\right) = \left\{\hskip.1em\raisebox{-3ex}{$\includegraphics{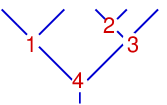}$}\hskip-.25em, \quad
	\raisebox{-3ex}{$\includegraphics{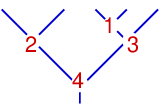}$}\hskip-.25em, \quad
	\raisebox{-3ex}{$\includegraphics{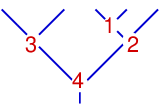}$}
\hskip.1em\right\}
\end{gather*}
\end{exam}

There are two right inverses to $\tau$ that will be useful later. Let $\Min(t)$ (respectively, $\Max(t)$) denote the unique $231$-avoiding ($132$-avoiding) permutation mapping to $t$ under $\tau$. Loday and Ronco show that $\tau^{-1}(t)$ is the interval $[\Min(t),\Max(t)]$ in the weak Bruhat order on the symmetric group \cite[Thm. 2.5]{LodRon:2002}, and that $\Min$ and $\Max$ are both order-preserving with respect to the Tamari order on $\y_n$.

\subsection{Bi-leveled trees and the multiplihedra}\label{sec: multitrees}

We next describe a family of \emph{bi-leveled trees} intermediate between the ordered and unordered ones. These trees arrange themselves as vertices of the \demph{multiplihedra}
$\m_\bb = \bigcup_n \m_n$, a family of polytopes introduced by Stasheff in 1970 \cite{Sta:1970} (though only proven to be polytopes much later \cite{For:2008}). Stasheff introduced this family to represent the fundamental structure of a weak map $f$ between weak structures, such as weak $n$-categories or $A_n$ spaces. The vertices of $\m_n$ correspond to associations of $n$ objects, pre- and post-application of $f$, e.g., $(f(a)f(b))f(c)$ and $f(a)f(bc)$. This leads to a natural description of $\m_n$ in terms of ``painted binary trees'' \cite{BoaVog:1973}, but we use here the description of Saneblidze and Umble \cite{SanUmb:2004}. 

A \demph{\mbox{bi-leveled} tree} is a pair $(t,C)$ with $t\in\y_n$ and $C\subseteq[n]$ designating some nodes of $t$ as lower than the others (indexing the nodes from left-to-right by $1,\ldots,n$). Viewing $t$ as a poset with root node maximal, $C$ is an increasing order ideal  in $t$ where the leftmost node is a minimal element. Graphically, $C$ indexes a collection of nodes of $t$ circled according to the rules: (i) the leftmost node is circled and has no circled children; (ii) if a node is circled, then its parent node is circled. 

Define a map $\beta$ from $\s_n$ to $\m_n$ as follows. Given a permutation $\sigma = \sigma_1\sigma_2\cdots\sigma_n$, first represent $\sigma$ as an ordered tree. Next, forget the ordering on the nodes, save for circling all nodes $\sigma_i$ with $\sigma_i\geq \sigma_1$. 

\begin{exam}\label{ex: lambda}
Consider again the permutations $1423$, $2413$, and $3412$ of Example \ref{ex: tau}. Viewed as ordered trees, their images under $\beta$ are distinct:
\begin{gather*}
	\left(\hskip.1em\raisebox{-3ex}{$\includegraphics{figures.FPSAC/p1423a_25}$}\hskip-.25em, \quad
	\raisebox{-3ex}{$\includegraphics{figures.FPSAC/p2413a_25}$}\hskip-.25em, \quad
	\raisebox{-3ex}{$\includegraphics{figures.FPSAC/p3412a_25}$}
\hskip.1em\right)
	\stackrel{\beta}{\longmapsto} 
	\left(\hskip.1em\raisebox{-3ex}{$\includegraphics{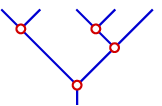}$}\hskip-.25em, \quad
	\raisebox{-3ex}{$\includegraphics{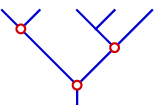}$}\hskip-.25em, \quad
	\raisebox{-3ex}{$\includegraphics{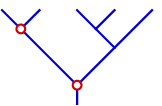}$}
\hskip.1em\right).
\end{gather*}
\end{exam}

\begin{figure}[!htb]
\centering
\scalebox{.90}{
\psset{unit=.5}
\begin{pspicture}(24,19.5)(.5,.5)
   \rput[c](12,10){\scalebox{.9}{\includegraphics{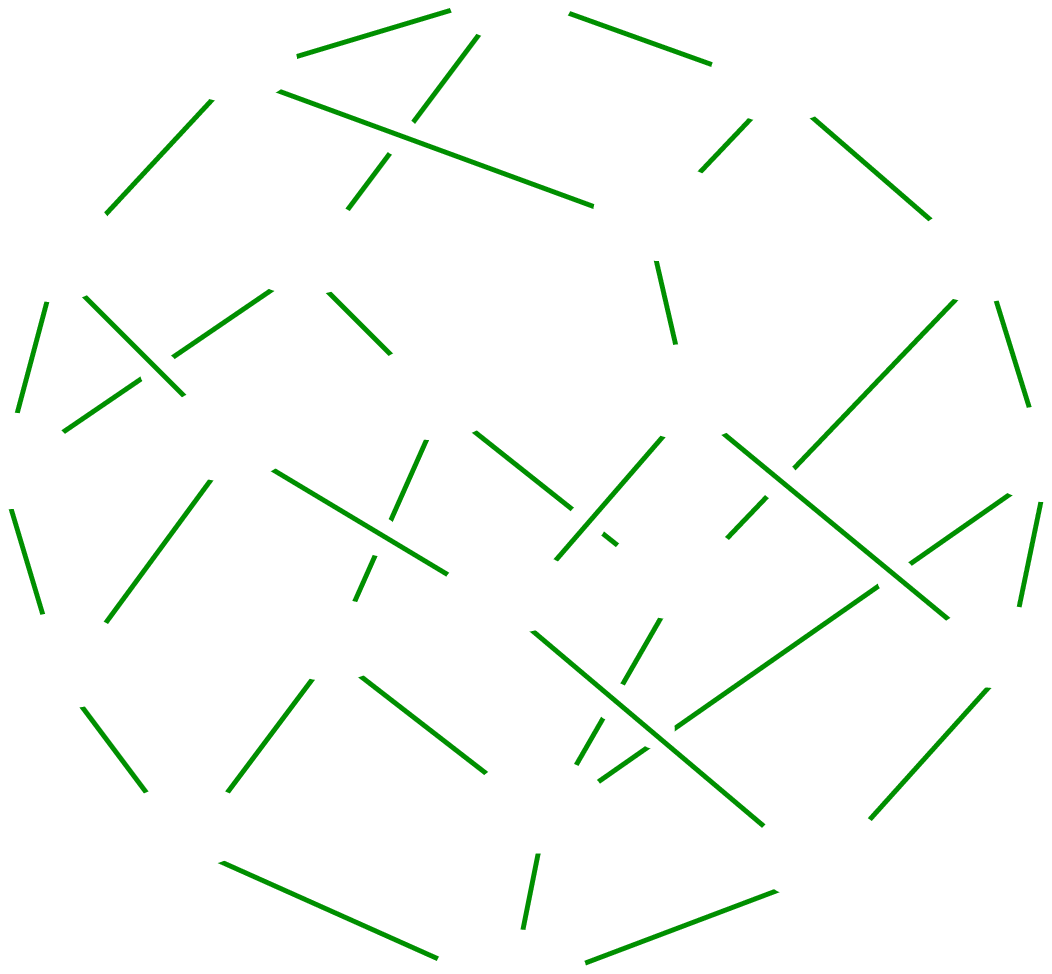}}}
 
   \rput[c](11.6,18.6){\includegraphics{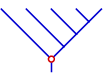}}

   \rput[c]( 6.6,17.5){\includegraphics{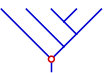}}
   \rput[c](16.6,17.2){\includegraphics{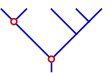}}

   \rput[c](3.5,14.0){\includegraphics{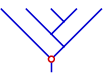}}
   \rput[c]( 7.8,14.1){\includegraphics{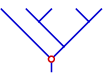}}
   \rput[c](14.55,14.75){\includegraphics{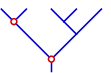}}
   \rput[c](20.25,14.0){\includegraphics{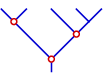}}

   \rput[c](2.3,10.15){\includegraphics{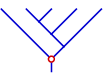}}
   \rput[c]( 6.6,10.65){\includegraphics{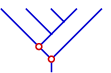}}
   \rput[c](10.5,11.45){\includegraphics{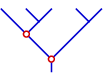}}
   \rput[c](15.0,11.5){\includegraphics{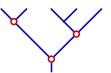}}
   \rput[c](21.55,10.45){\includegraphics{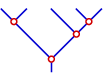}}

   \rput[c](3.5, 6.5){\includegraphics{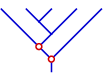}}
   \rput[c](8.5, 7){\includegraphics{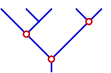}} 
   \rput[c](11.6,7.65){\includegraphics{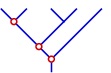}}
   \rput[c](14.75,8.1){\includegraphics{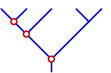}}
   \rput[c](20.8, 6.75){\includegraphics{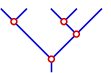}}

   \rput[c]( 5.75,3.5){\includegraphics{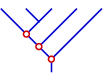}}
   \rput[c](12.5,3.8){\includegraphics{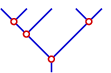}}
   \rput[c](17.5,3){\includegraphics{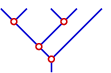}}

   \rput[c](11.75, .9){\includegraphics{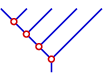}}
  \end{pspicture}}
\caption{The weak order on $\m_4$, the bi-leveled trees on $4$ nodes.}
\label{fig: M4}
\end{figure}

Denote by $\phi$ the map from bi-leveled trees to trees that forgets which nodes are circled. The map $\phi$ helps define a partial order on bi-leveled trees that extends the Tamari lattice on planar binary trees: say that the bi-leveled tree $s$ precedes the bi-leveled tree $t$ in the partial order if $\phi(s)\leq \phi(t)$ and the circled nodes satisfy $C_t\subseteq C_s$. We call this the weak order on bi-leveled trees. See Figure \ref{fig: M4} below for an example. 

\begin{figure}[!hbt]
\centering
\psset{unit=.335}
\begin{pspicture}(40,10.2)(0,.7)
%
\rput(20,5){%
$$
   \includegraphics[width=3.5cm]{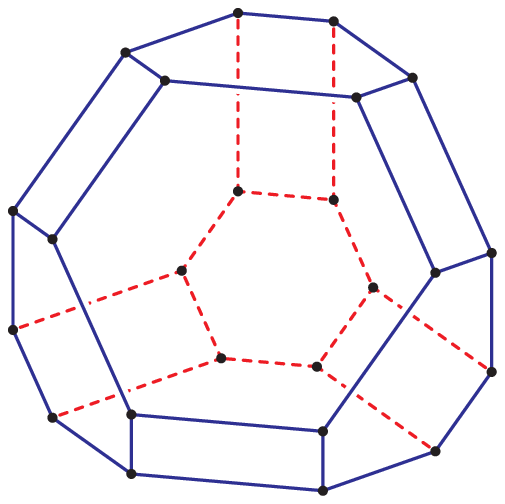}
\raisebox{11ex}{\large$\stackrel{\beta}{\relbar\joinrel\twoheadrightarrow}$\  \ }
   \includegraphics[width=3.5cm]{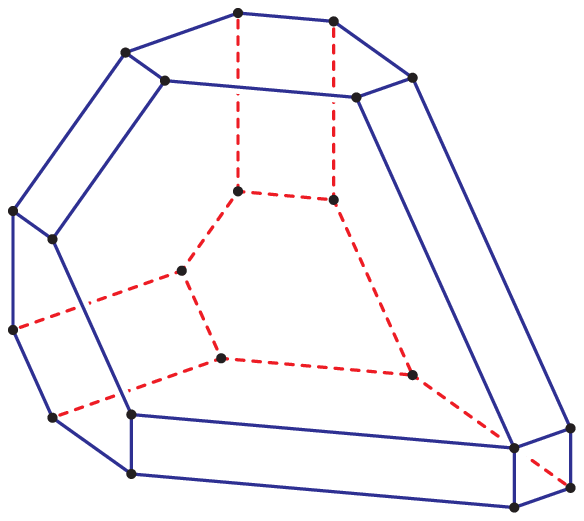}
\raisebox{11ex}{\large$\stackrel{\phi}{\relbar\joinrel\twoheadrightarrow}$\ \ }
   \includegraphics[width=3.5cm]{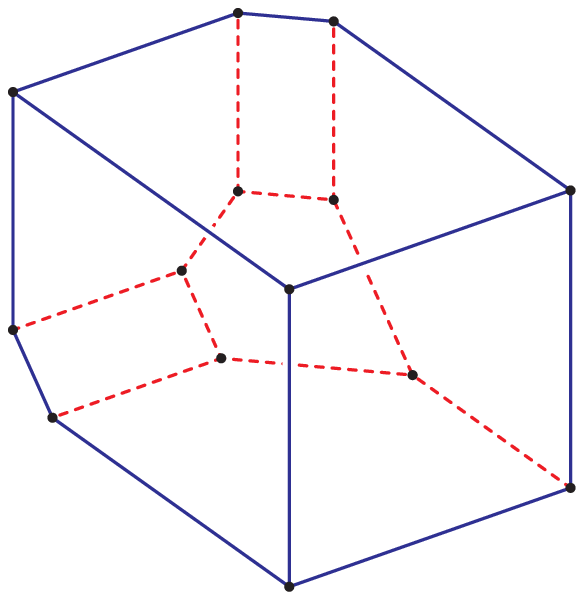}
$$}
{%
\pscircle*[linecolor=white,linewidth=.1](3.04,2.805){.21}
\pscircle*[linecolor=white,linewidth=.1](17.11,2.805){.21}
\pscircle*[linecolor=white,linewidth=.1](34.01,5.06){.21}
\pscircle[linecolor=green,linewidth=.08](3.04,2.805){.21}
\pscircle[linecolor=green,linewidth=.08](17.11,2.805){.21}
\pscircle[linecolor=green,linewidth=.08](34.01,5.06){.21}
}
\end{pspicture}
\caption{$\beta$ and $\phi$ extend to face (and poset) maps from the permutahedra to the associahedra. The distinguished vertices $1234$, $\beta(1234)$, and $\phi(\beta(1234))$ are indicated.}
\label{fig: polytopes}
\end{figure}

The equality $\phi\circ\beta = \tau$ is evident. Remarkably, this factorization $$\s_\bb \stackrel{\beta}{\relbar\joinrel\onto} \m_\bb \stackrel{\phi}{\relbar\joinrel\onto} \y_\bb$$ extends to the level of face maps between polytopes (see Figure \ref{fig: polytopes}). Our point of departure was the observation that it is also a factorization as poset maps. 

\subsection{Dimension enumeration}\label{sec: enum}

Fix a field $k$ of characteristic zero and let $\ssym$ denote the graded vector space $\bigoplus_{n\geq0} {\ssym}_{n}$ whose $n^{th}$ graded piece has the ``fundamental'' basis $\left\{F_\sigma \mid \sigma\right.$ an ordered tree in $\left.\Sym_n\right\}$. Define $\msym$ and $\ysym$ similarly, replacing $\s_n$ by $\m_n$ and $\y_n$, respectively. 
We follow convention and say that $\ssym_0$ and $\ysym_0$ are 1-dimensional. By contrast, we agree that $\msym_0 = \{0\}$. (See \cite{For:2} for categorical rationale; briefly, Stasheff's $\m_1$ is already $0$-dimensional, so $\m_0$ has no clear significance.)

In Section \ref{sec: msym}, we give these three vector spaces a variety of algebraic structures. Here we record some information about the dimensions of the graded pieces for later reference. 
\begin{align}
\label{eq: Hilb ssym}   \Hilb{q}{\ssym} &= \sum_{n\geq0} n!\,q^n = 1 + q + 2q^2 + 6q^3 + 24q^4 + 120q^5\cdots \\
\label{eq: Hilb msym}   \Hilb{q}{\msym} &=\sum_{n\geq1} A_n q^n =\hphantom{0+\,\,}q + 2q^2 + 6q^3 + 21q^4 + 80q^5 + \cdots \\
\label{eq: Hilb ysym}   \Hilb{q}{\ysym} &=\sum_{n\geq0} C_n q^n = 1 + q + 2q^2 + 5q^3 + 14q^4 + 42q^5 + \cdots 
\end{align}
Of course, $C_n$ is the $n^{th}$ Catalan number. The enumeration of bi-leveled trees is less familiar: the $n^{th}$ term satisfies $A_n = C_{n-1} + \sum_{k=1}^{n-1} A_i\,A_{n-i}$ \cite[A121988]{Slo:oeis}. 
A little generating function arithmetic can show that the quotient of \eqref{eq: Hilb msym} by \eqref{eq: Hilb ysym} expands as a power series with nonnegative coefficients,
\begin{gather}\label{eq: Hilb quotient}
\frac{\Hilb{q}{\msym}}{\Hilb{q}{\ysym}} = q+q^2+3 q^3+11 q^4+44 q^5+ \cdots .
\end{gather}
We will recover this with a little algebra in Section \ref{sec: coinvariants}. The positivity of the quotient of \eqref{eq: Hilb ssym} by \eqref{eq: Hilb ysym} is established by \cite[Theorem 7.2]{AguSot:2006}.

\section{The Hopf module $\msym$}\label{sec: msym}

Let $\btau$, $\bbeta$, and $\bphi$ be the maps between the vector spaces $\ssym$, $\msym$, and $\ysym$ induced by $\tau$, $\beta$, and $\phi$ on the fundamental bases. That is, for permutations $\sigma$ and bi-leveled trees $t$, we take
\begin{align*}
\btau(F_\sigma) &= F_{\tau(\sigma)}\,, & 
\bbeta(F_\sigma) &= F_{\beta(\sigma)}\,, & 
\bphi(F_t) &= F_{\phi(t)} \,. 
\end{align*}
Below, we recall the product and coproduct structures on the Hopf algebras $\ssym$ and $\ysym$. In \cite{MalReu:1995} and \cite{LodRon:1998}, these were defined in terms of the fundamental bases. Departing from these definitions, rich structural information was deduced about $\ssym$, $\ysym$, and the Hopf algebra map $\btau$ between them in \cite{AguSot:2005,AguSot:2006}. This information was revealed via a change of basis---from fundamental to ``monomial''---using M\"obius inversion. We take the same tack below with $\msym$ and meet with similar success.  

\subsection{The Hopf algebras $\ssym$ and $\ysym$}\label{sec: ssym}

Following \cite{AguSot:2006}, we define the product and coproduct structures on $\ssym$ and $\ysym$ in terms of \emph{$p$-splittings} and \emph{graftings} of trees. A \demph{$p$-splitting} of a tree $t$ with $n$ nodes is a forest (sequence) of $p+1$ trees with $n$ nodes in total. This sequence is obtained by choosing $p$ leaves of $t$ and splitting them (and all parent branchings) right down to the root. By way of example, consider the $3$-splitting below (where the third leaf is chosen twice and the fifth leaf is chosen once).
\[
\raisebox{13pt}{$t\,\,=$}\hskip-.1em  \includegraphics[width=0.8in]{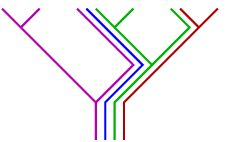}  
    \qquad\raisebox{13pt}{$\longrightarrow$}\qquad
   \raisebox{14pt}{$\Bigl($}
   \raisebox{7pt}{\begin{picture}(92,23)
    \put(3,0){\,\includegraphics[width=1.15in]{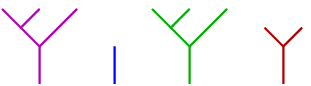}}
    \put(25,5){,} \put(42,5){,} \put(72,5){,}
   \end{picture}}
   \raisebox{14pt}{$\Bigr)$}\hskip-1.5em
    \qquad\raisebox{12pt}{$=\,\, (t_0,t_1,t_2,t_3)$ .}
\]
Denote a $p$-splitting of $t$ by $t \psplit (t_0,\ldots,t_{p})$. 
The \demph{grafting} of a forest $(t_0,t_1,\ldots,t_p)$ onto a tree with $p$ nodes is also best described in pictures; for the forest above and $s=\tau(213)$, the tree 
\begin{center}
  \includegraphics[width=0.8in]{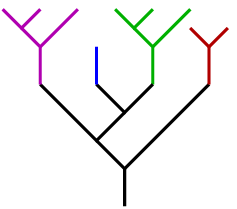}  
       \qquad\raisebox{15pt}{$=$}\qquad
  \raisebox{2pt}{\includegraphics[width=0.76in]{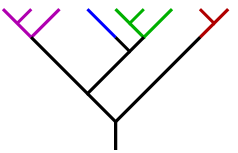}}
\end{center}
is the grafting of $(t_0,t_1,\ldots,t_p)$ onto $s$, denoted $(t_0,t_1,t_2,t_3)/s$.
Splittings and graftings of ordered trees are similarly defined. One remembers the labels originally assigned to the nodes of $t$ in a $p$-splitting, and if $t$ has $q$ nodes, then one increments the labels of $s$ by $q$ in a grafting $(t_0,t_1,t_2,t_3)/s$. See \cite{AguSot:2006} for details. 

\begin{defi}
Fix two ordered or ordinary trees $s$ and $t$ with $p$ and $q$ internal nodes, respectively. We define the product and coproduct by
\begin{equation}\label{eq: prod and coprod}
  F_t\bs\cdot F_s\ =\ 
    \sum_{t \psplit (t_0,t_1,\dotsc,t_p)} F_{(t_0,t_1,\dotsc,t_p)/s}
   \qquad\mbox{and}\qquad
  \Delta (F_t)\ =\ \sum_{t \psplit (t_0,t_1)} F_{t_0}\otimes F_{t_1}\,.
\end{equation}
(In the coproduct for ordered trees, the labels in $t_0$ and $t_1$ are reduced to be permutations of $|t_0|$ and $|t_1|$.)
\end{defi}

\subsection{Module and comodule structures}\label{sec: hopf module}

We next modify the structure maps in \eqref{eq: prod and coprod} to give $\msym$ the structure of (left) $\ssym$-module and (right) $\ysym$--Hopf module. Given a bi-leveled tree $b$, let $b\psplit(b_0,\ldots,b_p)$ represent any $p$-splitting of the underlying tree, together with a circling of all nodes in each $b_i$ that were originally circled in $b$.

\begin{defi}{\large\bf(}{\tiny\,}\emph{action} of $\ssym$ on $\msym${\tiny\,}{\large\bf)}{\ } 
For $w \in \s_\bb$ and $s\in \m_p$, write $b=\beta(w)$ and set
\begin{gather}\label{eq: action}
  F_w\bs\cdot F_s\ =\ 
    \sum_{b \psplit (b_0,b_1,\dotsc,b_p)} F_{(b_0,b_1,\dotsc,b_p)/s}
\end{gather}
where the circling rules in $(b_0,b_1,\dotsc,b_p)/s$ are as follows: every node originating in $s$ is circled whenever $|b_0|>0$, otherwise, every node originating in $b=\beta(w)$ is uncircled.
\end{defi}

This action may be combined with any section of $\beta$ to define a product on $\msym$. For example, 
 \begin{gather*}
F_{\includegraphics{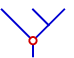}}  \bs\cdot F_{\includegraphics{figures.FPSAC/m21_25}} =
     {F_{\includegraphics{figures.FPSAC/m2143_2}}} + 
     {F_{\includegraphics{figures.FPSAC/m2413_2}}} + 
     {F_{\includegraphics{figures.FPSAC/m2431_2}}} +
     {F_{\includegraphics{figures.FPSAC/m4213_2}}} + 
     {F_{\includegraphics{figures.FPSAC/m4231_2}}} + 
     {F_{\includegraphics{figures.FPSAC/m4321_2}}}.
\end{gather*}
\begin{theo}
The action $\ssym\otimes \msym \to \msym$ and the product $\msym\otimes \msym \to \msym$ are associative. Moreover, putting $\msym_0 := k$, they make $\bbeta$ into an algebra map that factors $\btau$. 
\end{theo}

Unfortunately, no natural coalgebra structure exists on $\msym$ that makes $\bbeta$ into a Hopf algebra map.

\begin{defi}{\large\bf(}\emph{action}/\emph{coaction} of $\ysym$ on $\msym${\large\bf)}{\ } 
Given $b\in\m_\bb$, let $b\rsplit(b_0,\ldots,b_p)$ denote a $p$-splitting satisfying $|b_0|>0$. For $s\in \y_p$, set
\begin{equation}\label{eq: action and coaction}
  F_b\bs\cdot F_s\ =\ 
    \sum_{b \rsplit (b_0,b_1,\dotsc,b_p)} F_{(b_0,b_1,\dotsc,b_p)/s}\,
   \qquad\mbox{and}\qquad
  \brho (F_b)\ =\ \sum_{b \rsplit (b_0,b_1)} F_{b_0}\otimes F_{\phi(b_1)}\,,
\end{equation}
where in $(b_0,b_1,\dotsc,b_p)/s$ every node originating in $s$ is circled, and in $\phi(b_1)$ all circles are forgotten.
\end{defi}

\begin{exam}
In the fundamental bases of $\msym$ and $\ysym$, the action looks like
\begin{gather*}
F_{\includegraphics{figures.FPSAC/m21_25}}  \bs\cdot F_{\includegraphics{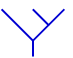}} =
     {F_{\includegraphics{figures.FPSAC/m2143_2}}} + {F_{\includegraphics{figures.FPSAC/m2413_2}}} + {F_{\includegraphics{figures.FPSAC/m2431_2}}} ,
\end{gather*}
while the coaction looks like 
\begin{align*}
\brho(F_{\includegraphics{figures.FPSAC/m3241_2}}) &=
     {F_{\includegraphics{figures.FPSAC/m3241_2}} \otimes 1} + 
     {F_{\includegraphics{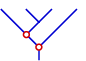}} \otimes  F_{\includegraphics{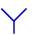}}} +   
     {F_{\includegraphics{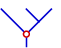}} \otimes  F_{\includegraphics{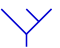}}} + 
     {F_{\includegraphics{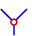}} \otimes  F_{\includegraphics{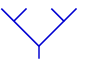}}} , \\
\brho(F_{\includegraphics{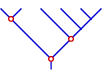}}) &=
     {F_{\includegraphics{figures.FPSAC/m35421_16}} \otimes 1} + 
     {F_{\includegraphics{figures.FPSAC/m2431_2}} \otimes  F_{\includegraphics{figures.FPSAC/1_2}}} + 
     {F_{\includegraphics{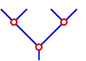}} \otimes  F_{\includegraphics{figures.FPSAC/21_2}}} + 
     {F_{\includegraphics{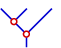}} \otimes  F_{\includegraphics{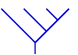}}}  +    
     {F_{\includegraphics{figures.FPSAC/m1_2}} \otimes  F_{\includegraphics{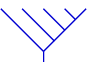}}} .
\end{align*}
\end{exam}

The significance of our definition of $\brho$ will be seen in Corollary \ref{thm: coinvariants}. Our next result requires only slight modifications to the original proof that $\ysym$ is a Hopf algebra (due to the restricted $p$-splittings). 

\begin{theo}\label{thm: Hopf module} The maps $\bs\cdot : \msym \otimes\ysym \to \msym$ and $\brho:\msym \to \msym \otimes \ysym$ are associative and coassociative, respectively. They give $\msym$ the structure of $\ysym$--Hopf module. That is, $\brho(F_b \bs\cdot F_s) = \brho(F_b) \bs\cdot \Delta(F_s)$. 
\end{theo}

\subsection{Main results}\label{sec: coinvariants}

We next introduce ``monomial bases'' for $\ssym$, $\msym$, and $\ysym$. Given $t\in \m_n$, define 
\[
	M_t = \sum_{t\leq t'} \mu(t,t') F_{t'},
\]
where $\mu(\,\cdot\,,\,\cdot\,)$ is the M\"obius function on the poset $\m_n$. Define the monomial bases of $\ssym$ and $\ysym$ similarly (see (13) and (17) in \cite{AguSot:2006}).
The coaction $\brho$ in this basis is particularly nice, but we need a bit more notation to describe it. 
Given $t\in \m_p$ and  $s\in \y_q$, let {$t\backslash s$} denote the bi-leveled tree on $p+q$ internal nodes formed by grafting the root of $s$ onto the rightmost leaf of $t$. 

\begin{theo}\label{thm: coaction} For a bi-leveled tree $t$, the coaction $\brho$ on $M_t$ is \, 
$\displaystyle
	\brho(M_t) = \sum_{t=t'\backslash s} M_{t'} \otimes M_s.
$
\end{theo}

\begin{exam}
Revisiting the trees in the previous example, the coaction in the monomial bases looks like
\begin{align*}
\brho(M_{\includegraphics{figures.FPSAC/m3241_2}}) &=
     {M_{\includegraphics{figures.FPSAC/m3241_2}} \otimes 1} + 
     {M_{\includegraphics{figures.FPSAC/m213_2}} \otimes M_{\includegraphics{figures.FPSAC/1_2}}}, \\%
\brho(M_{\includegraphics{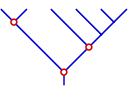}}) &=
     {M_{\includegraphics{figures.FPSAC/m35421_2}} \otimes 1} + 
     {M_{\includegraphics{figures.FPSAC/m2431_2}} \otimes  M_{\includegraphics{figures.FPSAC/1_2}}} + 
     {M_{\includegraphics{figures.FPSAC/m132_2}} \otimes  M_{\includegraphics{figures.FPSAC/21_2}}} .
\end{align*}
\end{exam}

Recall that the \demph{coinvariants} of a Hopf module $M$ over a Hopf algebra $H$ are defined by $M^{\mathrm{co}} = \left\{ m \in M \mid\right. $ $\left.\brho(m) = m\otimes 1 \right\}$. The fundamental theorem of Hopf modules provides that $M \simeq M^{\mathrm{co}} \otimes H$. The monomial basis of $\msym$ demonstrates this isomorphism explicitly. 

\begin{coro}\label{thm: coinvariants}A basis for the coinvariants in the Hopf module $\msym$ is given by $\bigl\{ M_t\bigr\}_{t\in \mathcal T}$, where $\mathcal T$ comprises the bi-leveled trees with no uncircled nodes on their right branches.
\end{coro}

This result explains the phenomenon observed in \eqref{eq: Hilb quotient}. It also parallels Corollary 5.3 of \cite{AguSot:2006} to an astonishing degree. There, the right-grafting idea above is defined for pairs of planar binary trees and used to describe the coproduct structure of $\ysym$ in its monomial basis. 

\section{Towards a proof of the main result}\label{sec: proof}

We follow the proof of \cite[Theorem 5.1]{AguSot:2006}, which uses properties of the monomial basis of $\ssym$ developed in \cite{AguSot:2005} to do the heavy lifting. In \cite{AguSot:2006}, the section $\y_\bb \stackrel{\Max}{\longrightarrow} \s_\bb$ of $\tau$ is shown to satisfy $\btau(M_{\Max(t)}) = M_t$ and $\btau(M_\sigma) = 0$ if $\sigma$ is not $132$-avoiding. This was proven using the following result about {Galois connections.}

\begin{theo}[{\cite[Thm. 1]{Rot:1964}}]\label{thm: Galois}
Suppose $P$ and $Q$ are two posets related by a \demph{Galois connection}, i.e., a pair of order-preserving maps $\varphi: P \to Q$ and $\gamma:Q \to P$ such that for any $v\in P$ and $t\in Q$,
$
	\varphi(v) \leq t \iff v \leq \gamma(t). 
$
Then the M\"obius functions $\mu_P$ and $\mu_Q$ are related by
\[
	\forall v\in P\hbox{ and }t\in Q, \quad
	\sum_{\newatop{w\in \varphi^{-1}(t),}{v\leq w}} \mu_P(v,w) \,=\, 
	\sum_{\newatop{s\in \gamma^{-1}(v),}{s\leq t}} \mu_Q(s,t).
\]
\end{theo}

There is a twist in our present situation. Specifically, no Galois connection exists between $\s_n$ and $\m_n$. On the other hand, we are not aiming for an order-preserving map $\iota:\m_\bb \into \s_\bb$ satisfying $\bbeta(M_{\iota(t)}) = M_t$ (indeed, no such map exists). Rather, we find that
\begin{gather}\label{eq: cell-sum}
	\bbeta\Bigg(\sum_{\sigma\in\beta^{-1}(t)} \! M_\sigma \Bigg) = M_t.
\end{gather}
This fact is the key ingredient in our proof of Theorem \ref{thm: coaction}. Its verification required modification of the notion of Galois connection---a relationship between posets that we call an \emph{\subGalois}(Section \ref{sec: subGalois}).

\subsection{Sections of the map $\beta:\s_\bb \to \m_\bb$}\label{sec: order-preserving}

Bi-leveled trees $t$ are in bijection with pairs $\{s,\mathbf s\}$, where $s$ is a planar binary tree, with $p$ nodes say, and $\mathbf s =(s_1,\ldots,s_{p})$ is a forest (sequence) of planar binary trees. In the bijection, $s$ comprises the circled nodes of $t$ and $s_i$ is the binary tree (of uncircled nodes) sitting above the $i^{\mathrm{th}}$ leaf of $s$. For example,
\[
	t = \hskip-.5em \raisebox{-4.5ex}{\includegraphics{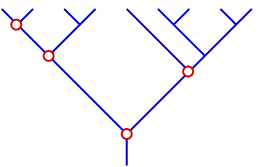}}
\,\,\longleftrightarrow\,\,\,\,
	\{s,\mathbf s\} = \left\{\raisebox{-2.75ex}{\includegraphics{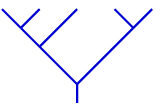}} ,\,\Big(
	\,\raisebox{-1ex}{\includegraphics{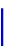}}\,,
	\raisebox{-1ex}{\includegraphics{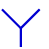}},
	\,\raisebox{-1ex}{\includegraphics{figures.FPSAC/0_3}}\,,
	\raisebox{-2.1ex}{\includegraphics{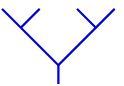}} \Big) \right\}.
\]

A natural choice for a section $\iota:\m_n \to \s_n$ would be to, say, build $\Min[s]$ and $\Min[s_i]$ for each $i$ and splice these permutations together in some way to build a word on the letters $\{1,2,\ldots n\}$. Let  $\minmin(t)$ denote the choice giving $s_1$ smaller letters than $s_2$, $s_2$ smaller letters than $s_3$, \dots, $s_{p-1}$ smaller letters than $s_p$, and $s_p$ smaller letters than $s$:
\[
	\minmin\hskip-.2ex \Biggl(\,\raisebox{-4ex}{\includegraphics{figures.FPSAC/fromUtoVexam05_25}}\,\Biggr) = \minmin\hskip-.2ex \Biggl(\,\raisebox{-4ex}{\includegraphics{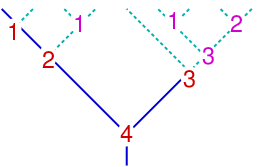}}\,\Biggr) = 
	\raisebox{-4ex}{\includegraphics{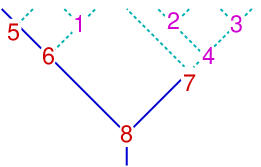}}  \!\!=\,\, 56187243.
\]
This choice does not induce a poset map. The similarly defined $\maxmax$ also fails (chosing maximal permutations representing $\bs s$ and $s$), but $\maxmin$ has the properties we need: 
\[
	\maxmin\hskip-.2ex  \Biggl(\,\raisebox{-4ex}{\includegraphics{figures.FPSAC/fromUtoVexam05_25}}\,\Biggr) = 
	\raisebox{-4ex}{\includegraphics{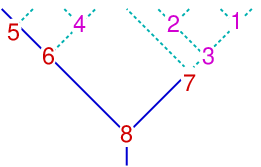}}  \!\!=\,\, 56487231.
\]

We define this map carefully. Given  $t\in \y_n$ and any subset $S\subseteq \N$ of cardinality $n$, write $\Min[]_S(t)$ for the image of $\Min[t]$ under the unique order-preserving map from $[n]$ to $S$; define $\Max[]_S(t)$ similarly. 

\begin{defi}[The section $\maxmin$]\label{def: iota}
Let $t\leftrightarrow\{s,\mathbf{s}\}$ be a bi-leveled tree on $n$ nodes with $p$ circled nodes. Write $ u= u_1\cdots u_p = \Min[]_{[a,b]}(s)$ for $[a,b]=\{n-p+1,\ldots,n\}$ and write $v^i=\Max[]_{[a_i,b_i]}(s_i)$ ($1\leq i \leq p$), where the intervals $[a_i,b_i]$ are defined recursively as follows:
\begin{align*}
	a_p &= 1 \quad\hbox{and}\quad b_p = a_p+|s_p|-1, \\
	a_i &= 1+ \max \bigcup_{j>i} S_j \quad\hbox{and}\quad  b_i = a_i + |s_i| - 1.
\end{align*}
Finally, define $\maxmin(t)$ by the concatenation 
$
	\maxmin(t) =  u_1 v^1 u_2 v^2\cdots  u_p v^p. 
$
\end{defi}

\begin{rema}Alternatively, 
$\maxmin(t)$ is the unique $w\in \beta^{-1}(t)$ avoiding the pinned patterns 
$\underline{0}231,$ $\underline{3}021,$ and $\underline{2}031,$
where the underlined letter is the first letter in $w$. The first two patterns fix the embeddings of $s_i$ ($0\leq i \leq p$), the last one makes the letters in $s_i$ larger than those in $s_{i+1}$ ($1\leq i < p$).
\end{rema}

The important properties of $\maxmin(t)$ are as follows.

\begin{prop}\label{thm: beta-maxmin} The section $\iota: \m_n \to \s_n$ given by $\iota(t) = \maxmin(t)$ is an embedding of posets. The map $\beta: \s_n \to \m_n$ satisfies $\beta(\iota(t))=t$ \,for all \,$t\in \m_n$ and $\beta^{-1}(t) \subseteq \s_n$ is the interval $[\minmin(t),\maxmax(t)]$.
\end{prop}

\subsection{Interval retracts}\label{sec: subGalois}

Let $\varphi:P \to Q$ and $\gamma: Q \to P$ be two order-preserving maps between a lattice $P$ and a poset $Q$. If
\begin{gather*}
	\forall t\in Q \qquad
	\varphi(\gamma(t)) = t \quad\hbox{and}\quad
	\varphi^{-1}(t) \ \hbox{ is an interval}, \ 
\end{gather*}
then we say that $\varphi$ and $\gamma$ demonstrate $P$ as an \demph{\subGalois}of $Q$.

\begin{theo}\label{thm: subGalois}
If $P$ and $Q$ are two posets related by an {\subGalois} $(\varphi,\gamma)$, then the M\"obius functions $\mu_P$ and $\mu_Q$ are related by
\[
	\forall s<t\in Q \quad
	\sum_{\newatop{v\in \varphi^{-1}(s)}{w\in \varphi^{-1}(t)}} \! \mu_P(v,w) \,=\, 
	\mu_Q(s,t).
\]
\end{theo}

The proof of Theorem \ref{thm: subGalois} exploits Hall's formula for M\"obius functions. An immediate consequence is a version of \eqref{eq: cell-sum} for any $P$ and $Q$ related by an \subGalois. Verifying that $(\beta,\maxmin)$ is an \subGalois between $\s_n$ and $\m_n$ (Proposition \ref{thm: beta-maxmin}) amounts to basic combinatorics of the weak order on $\s_n$.

\section{More families of binary trees and their polytopes}\label{sec: polytopes}
We have so far ignored the algebra $\qsym$ of quasisymmetric functions advertised in the introduction. A basis for its $n^{th}$ graded piece is naturally indexed by compositions of $n$, but may also be indexed by trees as follows. To a composition $(a_1,a_2,\ldots)$, say $(3,2,1,4)$, we associate a sequence of \demph{right-combs} 
\[
(\raisebox{-0.75ex}{\includegraphics{figures.FPSAC/21_2}}\!,\ \raisebox{-0.70ex}{\includegraphics{figures.FPSAC/1_2}},\  \raisebox{-0.65ex}{\includegraphics{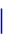}}{\scriptscriptstyle\,}, \raisebox{-0.70ex}{\includegraphics{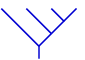}}),
\]
i.e., trees with $a_i$ leaves and all internal leaves rooted to the rightmost branch and left-pointing. These may be hung on another tree, a left-comb with right-pointing leaves, to establish a bijection between compositions of $n$ and ``combs of combs'' with $n$ internal nodes:
\begin{center}
\psset{unit=.5}
\begin{pspicture}(28,3.60)(.50,.25)
\rput[c](14,1.4){$
(23) \leftrightarrow \raisebox{-0.55ex}{\includegraphics{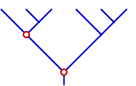}} \qquad 
(32) \leftrightarrow \raisebox{-0.55ex}{\includegraphics{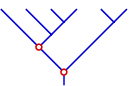}} \qquad 
(1214) \leftrightarrow \raisebox{-3.3ex}{\includegraphics{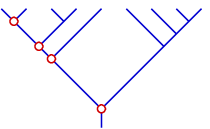}}  \hskip-0.5em.
$}
\rput[c](6.99,3.10){\footnotesize$\underline{\, 2\, }$}
\rput[c](8.21,3.10){\footnotesize$\underline{\,\ 3\ \,}$}
\rput[c](13.95,3.10){\footnotesize$\underline{\,\ 3\ \,}$}
\rput[c](15.22,3.10){\footnotesize$\underline{\, 2\, }$}
\rput[c](20.92,3.10){\footnotesize$\underline{1}$}
\rput[c](21.68,3.10){\footnotesize$\underline{\,2\,}$}
\rput[c](22.42,3.10){\footnotesize$\underline{1}$}
\rput[c](23.72,3.10){\footnotesize$\underline{\;\ \ 4\ \ \;}$}
\end{pspicture}
\end{center}
To see how $\qsym$ and the hypercubes fit into the picture, we briefly revisit the map $\beta$ of Section \ref{sec: multitrees}. 

We identified bi-leveled trees with pairs $\{s,\mathbf{s}\}$, where $s$ is the tree of circled nodes and $\mathbf s =(s_0,\ldots,s_{p})$ is a forest of trees (the uncircled nodes). Under this identification, $\beta$ may be viewed as a pair of maps $(\tau,\tau)$---with the first factor $\tau$ making a (planar binary) trees out of the nodes greater than or equal to $\sigma_1$, and the second factor $\tau$ making trees out of the smaller nodes:
\begin{gather*}
	\raisebox{-2ex}{$\includegraphics{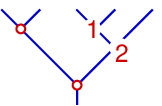}$} \ \ 
	\stackrel{(\tau,\id)}{\longleftarrow\joinrel\relbar\joinrel\relbar} \ \ 
	\raisebox{-2ex}{$\includegraphics{figures.FPSAC/p3412a_25}$} \ \ 
	\stackrel{(\id,\tau)}{\relbar\joinrel\relbar\joinrel\longrightarrow} \ \ 
	\raisebox{-2ex}{$\includegraphics{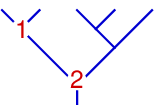}$}
	\hskip-0.5em.
\end{gather*}
See also Figure \ref{fig: tree-maps}.
Two more fundamental maps are $\gamma_l$ and $\gamma_r$, taking trees to (left- or right-) combs, e.g.,

\[
	\raisebox{-0.55ex}{\includegraphics{figures.FPSAC/231_2}} \ 
	\stackrel{\gamma_l}{\longrightarrow} \ 
	\raisebox{-0.55ex}{\includegraphics{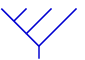}} \qquad \hbox{and}\qquad
	\raisebox{-1.3ex}{\includegraphics{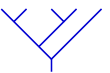}} \ 
	\stackrel{\gamma_r}{\longrightarrow} \ 
	\raisebox{-1.3ex}{\includegraphics{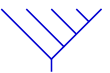}} 
	\hskip-0.25em.
\]
Figure \ref{fig: tree-maps} displays several combinations of the maps $\tau,$ $\gamma_l,$ and $\gamma_r$. The algebra $\qsym$ corresponds to the terminal object there---the set denoted 
$\bigl\{\frac{\cdots\hbox{\footnotesize\  combs }\cdots}{\hbox{\footnotesize comb}}\bigr\}$.

The new binary tree--like structures appearing in the factorization of $\ssym \onto \qsym$ (i.e., those trees not appearing on the central, vertical axis of Figure \ref{fig: tree-maps}) will be studied in upcoming papers. 
It is no surprise that $\ssym \onto \qsym$ factors through so many intermediate structures. What is remarkable, and what our binary tree point-of-view reveals, is that each family of trees in Figure \ref{fig: tree-maps} can be arranged into a family of polytopes. See Figure \ref{fig: polytope-maps}.  The (Hopf) algebraic and geometric implications of this phenomenon will also be addressed in future work. 

\begin{figure}[!hbt]
{\parindent-3em
\include{figures.FPSAC/hopftrees}}
\caption{A commuting diagram of tree-like objects. The spaces $\ssym$, $\msym$, $\ysym$ and $\qsym$ appear, top to bottom, along the center. 
The unlabeled dashed line represents the usual map from $\ysym$ to $\qsym$ (see \cite{LodRon:1998}, Section 4.4). It is incompatible with the given map $(\gamma,\gamma){\,\,:\,\,}\msym \to \qsym$. 
{\color{white} $!_{!_{!_{!_{!_{!_{!_{!}}}}}}}$}
\newline
We explore the Hopf module structures of objects mapping to $\ysym$ and $\qsym$ in future work. At least some of these will be full-fledged Hopf algebras (e.g., note that there is a bijection of sets between $\bigl\{\frac{\cdots\hbox{\footnotesize\  trees }\cdots}{\hbox{\footnotesize comb}}\bigr\}$ and $\bigl\{{\hbox{\footnotesize trees}}\bigr\}$, the latter indexing the Hopf algebra $\ysym$). 
}
\label{fig: tree-maps}
\end{figure}

\begin{figure}[!hbt]
\vskip3ex

\hskip1.5ex$\scalebox{.75}{\includegraphics{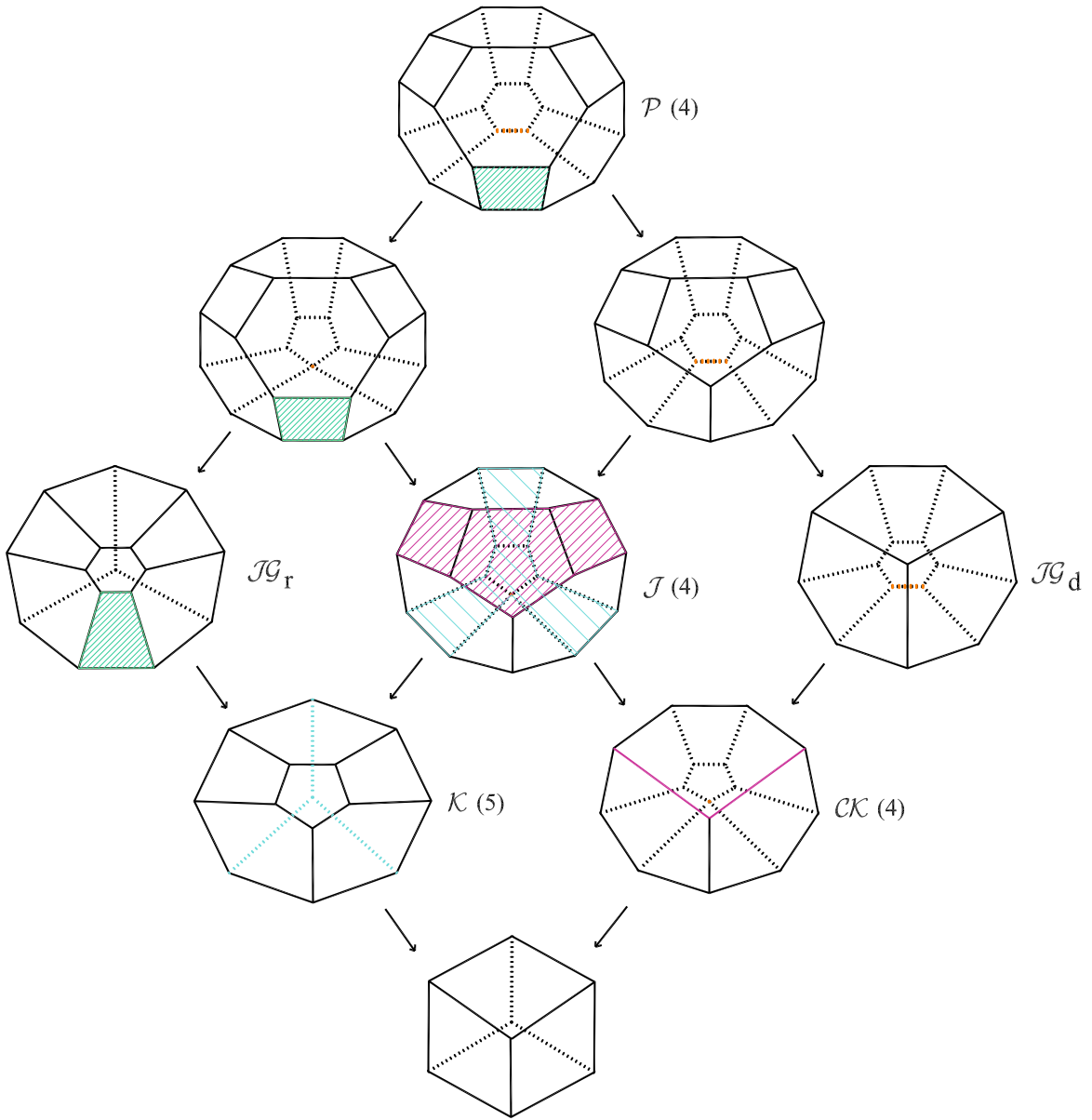}}$

\vskip1ex

\caption{A commuting diagram of polytopes based on tree-like objects with $4$ nodes, corresponding position-wise to Figure \ref{fig: tree-maps} (image of $\phi$ is suppressed). 
Notation is taken from \cite{DevFor:1}: $\mathcal{P}(4)$ is the permutohedron, 
$\mathcal{J}(4)$ is the multiplihedron, 
$\mathcal{K}(5)$ is the associahedron, and
$\mathcal{CK}(4)$ is the composihedron.
$\mathcal{JG}_d$ is the domain quotient of the permutohedron and 
$\mathcal{JG}_r$ is its range quotient. 
{\color{white} $!_{!_{!_{!_{!_{!_{!_{!}}}}}}}$}
\newline
The cellular projections shown include neither the Tonks projection nor the Loday Ronco projection from the associahedron to the hypercube. However, the map from the multiplihedron to the cube passing through the associahedron appears in \cite{BoaVog:1973}.
}
\label{fig: polytope-maps}
\end{figure}



\small

\providecommand{\bysame}{\leavevmode\hbox to3em{\hrulefill}\thinspace}
\providecommand{\MR}{\relax\ifhmode\unskip\space\fi MR }
\providecommand{\MRhref}[2]{%
  \href{http://www.ams.org/mathscinet-getitem?mr=#1}{#2}
}
\providecommand{\href}[2]{#2}

\end{document}

%% file: figures.FPSAC/hopftrees.tex
$
\xymatrix@R=12pt @C=0pt @M=2pt{
& & *++{\begin{array}{c} \bigl\{\hbox{\footnotesize permutations}\bigr\} \\[2ex] 
 \raisebox{-.5\height}{\scalebox{1.1}{\includegraphics{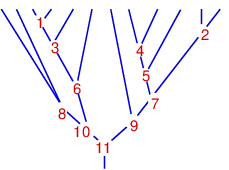}}} 
\end{array}} 
	\ar@{->}[ddl]_{(\tau,\id)} \ar@{->}[ddr]^{(\id,\tau)}  \ar@{-->}[dddd]^{\displaystyle\beta}  \\ \\
& *++{\begin{array}{c}  \bigl\{\frac{\cdots\hbox{\footnotesize\  perms }\cdots}{\hbox{\footnotesize tree}} \bigr\} \\[2ex] 
\raisebox{-.5\height}{\scalebox{1.0}{\includegraphics{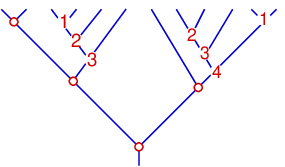}}}
\end{array}}
	\ar@{->}[ddl]_{(\gamma_r,\id)} \ar@{->}[ddr]^{(\id,\tau)}   && 
*++{\begin{array}{c}  \bigl\{\frac{\cdots\hbox{\footnotesize\  trees }\cdots}{\hbox{\footnotesize perm}} \bigr\} \\[2ex] 
\raisebox{-.5\height}{\scalebox{1.0}{\includegraphics{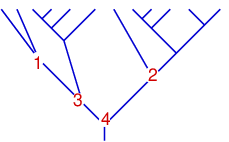}}}
\end{array}}
	\ar@{->}[ddl]_{(\tau,\id)} \ar@{->}[ddr]^{(\id,\gamma_l)}   \\ \\
*++{\begin{array}{c}  \bigl\{\frac{\cdots\hbox{\footnotesize\  perms }\cdots}{\hbox{\footnotesize comb}}\bigr\}  \\[2ex] 
\raisebox{-.5\height}{\scalebox{1.0}{\includegraphics{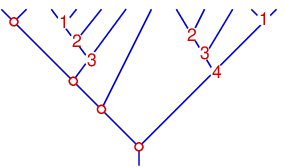}}}
\end{array}} 
	\ar@{->}[ddr]^{(\id,\tau)}   && 
*++{\begin{array}{c} \bigl\{\frac{\cdots\hbox{\footnotesize\  trees }\cdots}{\hbox{\footnotesize tree}} \bigr\} \\[2ex]  
\raisebox{-.5\height}{\scalebox{1.0}{\includegraphics{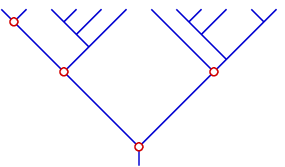}}}  \end{array}} 
	\ar@{->}[ddl]_{(\gamma_r,\id)} \ar@{->}[ddr]^{(\id,\gamma_l)} \ar@{-->}[dd]^{\raisebox{2.5ex}{$\displaystyle\phi$}} && 
*++{\begin{array}{c}  \bigl\{ \frac{\cdots\hbox{\footnotesize\  combs }\cdots}{\hbox{\footnotesize perm}}\bigr\} \\[2ex] 
\raisebox{-.5\height}{\scalebox{1.0}{\includegraphics{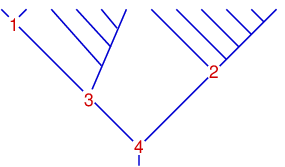}}}
\end{array}}  
	\ar@{->}[ddl]_{(\tau,\id)}    \\ \\
& *++{\begin{array}{c} \bigl\{\frac{\cdots\hbox{\footnotesize\  trees }\cdots}{\hbox{\footnotesize comb}}\bigr\} \\[2ex]   
\raisebox{-.5\height}{\scalebox{1.0}{\includegraphics{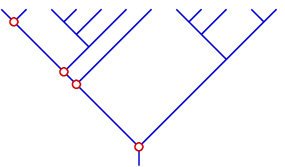}}}   \end{array}} 
	 \ar@{->}[ddr]^{(\id,\gamma_l)}  &
*++{\begin{array}{c} \bigl\{\hbox{\footnotesize trees}\bigr\} \\[2ex]  
\raisebox{-.5\height}{\scalebox{0.9}{\includegraphics{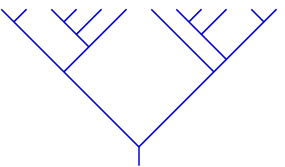}}}  \end{array}}
	 \ar@{-->}[dd] & 
*++{\begin{array}{c} \bigl\{\frac{\cdots\hbox{\footnotesize\  combs }\cdots}{\hbox{\footnotesize tree}} \bigr\} \\[2ex]  
\raisebox{-.5\height}{\scalebox{1.0}{\includegraphics{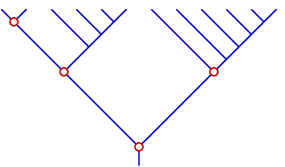}}}   \end{array}}    
	\ar@{->}[ddl]_{(\gamma_r,\id)}   \\ \\
& & *++{\begin{array}{c} \bigl\{\frac{\cdots\hbox{\footnotesize\  combs }\cdots}{\hbox{\footnotesize comb}}\bigr\} \\[2ex]  
\raisebox{-.5\height}{\scalebox{1.0}{\includegraphics{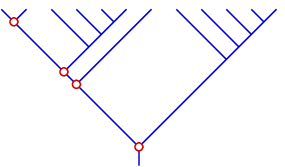}}}  \end{array}} &     
}
$